\documentclass{article}
\font\Bb=msbm10
\def\CC{\hbox{\Bb C}}

\def\RR{\hbox{\Bb R}}
\def\PP{\hbox{\Bb P}}

\usepackage{amssymb,latexsym,amsmath}

\newtheorem{theorem}{Theorem}[section]

\newtheorem{proposition}[theorem]{Proposition}

\title{When do linear combinations of orthogonal polynomials 
yield new sequences of orthogonal polynomials?}

\author{Manuel Alfaro \thanks{Partially supported by
Ministerio de Educaci\'on y Ciencia (MEC) of Spain under Grant MTM 2006-13000-C03-03 and
Diputaci\'on General de Arag\'on (DGA) project E-64.}
\\ Dpto. de Matem\'{a}ticas and IUMA. Universidad de Zaragoza (Spain)
\and  Francisco Marcell\'{a}n \thanks{Partially supported by MEC of Spain under Grant MTM
2006-13000-C03-02 and INTAS Research Network NeCCA 03-51-6637.} \\Dpto. de
Matem\'{a}ticas.  Universidad Carlos III de Madrid (Spain) \\
\and  Ana Pe\~{n}a \thanks{Partially supported by MEC of Spain under Grants MTM 2004-03036
and MTM 2006-13000-C03-03 and DGA project E-64.} \\Dpto. de Matem\'{a}ticas.
Universidad de Zaragoza (Spain) 
\and
M. Luisa Rezola$^{*}$ \\Dpto. de Matem\'{a}ticas and IUMA. Universidad de Zaragoza
(Spain)}

\date{}

\begin{document}

\maketitle
\vspace*{-12pt}
Dedicated to Professor Jes\'us S. Dehesa on the occasion of his 60th birthday

\begin{abstract}
Given $\{ P_n \}_{n \ge 0}$ a sequence of monic orthogonal
polynomials, we analyze their linear combinations with constant
coefficients and fixed length, i.e.,

$Q_n(x)=P_n(x)+a_1P_{n-1}(x)+\dots +a_kP_{n-k},\, a_k\not=0, \,n > k .$

\noindent Necessary and sufficient conditions are given for the orthogonality of
the sequence $\{ Q_n \}_{n \ge 0}$ as well as an interesting
interpretation in terms of the Jacobi matrices associated with $\{
P_n \}_{n \ge 0}$  and $\{ Q_n \}_{n \ge 0}$.

Moreover, in the case $k=2$, we characterize the families $\{ P_n
\}_{n \ge 0}$ such that the corresponding polynomials $\{ Q_n \}_{n \ge 0}$ are also
orthogonal.
\end{abstract}






\section{Introduction and basic definitions}
\setcounter{equation}{0} Given a linear functional $u$ on the linear
space $\PP$ of polynomials with real coefficients, a sequence of
monic polynomials $\{ P_n \}_{n \ge 0}$ with deg $P_n = n$ is said
to be orthogonal with respect to $u$ if $\langle u, P_n P_m \rangle = 0 $ for every $n
\not=m$ and $\langle u, P_n^2 \rangle \not= 0 $ for every $n = 0, 1, \dots \, .$

A linear functional $u$ is said to be quasi--definite (respectively
positive definite) if the leading principal submatrices $H_n$ of the Hankel matrix $H=(u_{i+j})_{i,j\ge0}$
associated with $u$, where $u_k = \langle u, x^k \rangle \,, k \ge 0 \,,$ are nonsingular
(respectively positive definite) for every $n \ge 0 $ (see \cite{Ch}).

A very well known result (Favard's theorem, see \cite{Ch} for instance) gives
a characterization of a quasi--definite (respectively positive
definite) linear functional in terms of the three--term recurrence
relation that the sequence $\{ P_n \}_{n \ge 0}$ satisfies, i.e.
\begin{align} \label{recurrenciaPn}
xP_n(x)&=P_{n+1}(x) + \beta_n P_n(x) + \gamma_n P_{n-1}(x), \quad
\\ P_0(x)&=1, \quad P_1(x)=x-\beta_0,\nonumber
\end{align}
whith $\gamma_n \not=0$ (respectively $\gamma_n
> 0)$.

In particular, if $u$ is a positive definite linear functional then
there exists a positive Borel measure $\mu$ supported on an infinite
subset of $\RR$ such that $\langle u, q \rangle = \int_{\RR} q \,
d\mu$ for every $q \in \PP \,.$
In such a situation, the zeros of $P_n$ are real, simple, and they are located in the
convex hull of the support of the measure $\mu$. Furthermore, the
zeros of $P_{n-1}$ interlace with those of $P_n$. 
Actually, this is a relevant fact in numerical quadrature, i.e. in the
discrete representation
\begin{equation} \label{cuadratura}
\int_{\RR} q \, d\mu \sim \sum_{k=1}^n \lambda_k q(c_k) \,, \quad
q \in \PP \,.
\end{equation}
If we choose as $(c_k)_{k=1}^n$ the zeros of $P_n$ then
(\ref{cuadratura}) is exact for every polynomial of degree at most
$2n - 1$ and, as a consequence of the interlacing property aforementioned, the
Christoffel-Cotes numbers $(\lambda_k)_{k=1}^n$ are positive real numbers.

In general, given the pair $(q, \mu)$ with $\displaystyle q(x) =
\prod_{k=1}^n (x - c_k)$ and letting $\lambda = (\lambda_1, \dots ,
\lambda_n)$ where $\displaystyle \lambda_k = \int_{\RR}
\frac{q(x)}{q'(c_k)(x - c_k)}\, d\mu(x)$,  $1 \le k \le n$, there
exists an integer number $d(q, \mu)$ with $n - 1 \le d(q, \mu)
\le 2n - 1 \,,$ so that (\ref{cuadratura}) is exact for the polynomials of
degree $\le d(q, \mu)$ but not for all polynomials of degree $d(q,
\mu) + 1 \,.$ The number $d(q, \mu)$ is said to be the degree of
precision of $(q, \mu).$

Shohat, in \cite{Sh}, proved that $(q, \mu)$ has degree of
precision $2n - 1 - k$ if and only if $q = P_n + a_1 P_{n-1} +
\dots  + a_k P_{n-k}$ where $a_k \not= 0$ and $\{ P_n \}_{n \ge 0}$ is the
sequence of monic polynomials orthogonal with respect to the measure
$\mu$.

Moreover, when  $\rm{supp} \, \mu = (-1, 1)$, Peherstorfer addresses in \cite{Peh}
sufficient conditions on the real numbers $\{ a_j \}_{j=1}^k$ under which 
the polynomial $q = P_n + a_1 P_{n-1} + \dots + a_k P_{n-k}$ has $n$ simple
zeros in $(-1, 1)$ and whose Christoffel-Cotes numbers are positive.

In \cite{Sh} a discussion about the zeros of the polynomial $q =
P_n + a_1 P_{n-1}$ is given in terms of $\rm{sign} \, a_1$: they are
real and simple and at most one of them lies outside
$\rm{supp} \, \mu$. Moreover, the zeros of the polynomial $q = P_n + a_1 P_{n-1} +
a_2 P_{n-2}$ are studied. If $a_2 < 0$, all the zeros are real and
simple and at most two of them do not belong to the $\rm{supp} \, \mu$.
In addition, in \cite{BDR-Z} it is proved that if $a_2 < 0$ then the
zeros of $P_{n-1}$ interlace with the zeros of $q$. The position of
the least and greatest zero of $q$ in terms of the least and
greatest zero of $P_n$ is also analyzed.

In \cite{BD} the positivity of Christoffel-Cotes numbers and the distribution of
zeros of linear combinations $R = P_m + \dots  + a_s P_s$ where $a_s
\not= 0 \,, 1 \le s \le m \le n$ and $m \le d(q, \mu)$ is analyzed.
Here $\displaystyle q(x) = \prod_{k=1}^n (x - c_k)$ with $c_1 <
\dots < c_n \,.$
If all the Christoffel-Cotes numbers are positive, then either $R$ is a non--zero scalar
multiple of $q$ or at least $N$ of the intervals $(c_k, c_{k+1})$ contain a zero of $R$
where
$N = \min \{s, d(q, \mu) + 1 -m \} \ge 1 \,.$

Grinshpun, in \cite{G}, studied the orthogonality of
special linear combinations of polynomials orthogonal with respect
to a weight function supported on an interval of the real line.
Such families of orthogonal polynomials come up in some extremal
problems of Zolotarev--Markov type as well as in problems of least
deviating from zero. He proved that the Bernstein--Szeg\H{o}
polynomials can be represented as a linear combination of the Chebyshev polynomials of
the same kind. Nevertheless, the special feature of this representation is that the
coefficients do not depend on $n$. The relevant question is if this property
characterizes Bernstein--Szeg\H{o} polynomials. Theorem 3.1 in
\cite{G} gives a positive answer in the sense that
Bernstein--Szeg\H{o} polynomials and just them can be represented as a linear
combination of Chebyshev polynomials with constant coefficients
independent of $n$ and fixed length. In other words, $\{ Q_n \}_{n\ge0}$ with
$Q_n = P_n + a_1 P_{n-1} + \dots + a_k P_{n-k}\,, n>k \,,$ where $\{ P_n \}_{n\ge0}$
is the Chebyshev sequence of j--th kind $(j=1,2,3,4)$ and $a_k \not= 0$, is a
sequence of orthogonal polynomials with respect to a weight
$\widetilde{\omega}$ if and only if $\displaystyle \widetilde{\omega}(x) =
\frac{\mu_j(x)}{h_k(x)} \,,$ where $h_k$ is a polynomial
of degree $k$ positive on $(-1, 1)$ and $\mu_j$ is the Chebyshev weight of
j--th kind, $(j = 1, 2, 3, 4)$.

\medskip
The aim of this work is to analyze linear combinations with constant
coefficients $Q_n = P_n + a_1 P_{n-1} + \dots + a_k P_{n-k}\,, n>k \,,$ of a sequence of
orthogonal polynomials $\{ P_n \}_{n\ge0}$. In Section 2 we find necessary and
sufficient conditions so that the sequence  $\{ Q_n \}_{n\ge0}$
is orthogonal with respect to a linear functional $v$.
Moreover, we discuss the matrix representation for the multiplication operator
in terms of the bases $\{ P_n \}_{n\ge0}$ and $\{ Q_n \}_{n\ge0}$
respectively. Such a matrix is a monic tridiagonal (Jacobi) matrix. We
prove that the leading principal submatrix associated with $\{ Q_n
\}_{n\ge0}$ is similar to a rank--one perturbation of the leading
principal submatrix associated with $\{ P_n \}_{n\ge0}.$  Also, we
give a simple algorithm to compute the polynomial $h_k$ of degree $k$ appearing in the
relation between the two functionals, $u=h_k v$.

In Section 3, the case $k=2$ is addressed, describing all the families $\{P_n \}_{n\ge0}$ 
orthogonal with respect to a linear functional such that the
corresponding $\{ Q_n \}_{n\ge0}$ is also orthogonal, obtaining explicit expressions for
the recurrence parameters
$\{ \beta_n \}_{n\ge0}$ and $\{ \gamma_n \}_{n\ge1}$ of the sequence $\{P_n \}_{n\ge0}$.
Finally, in Section 4 we present some remarks and examples of such
sequences
$\{ P_n \}_{n\ge0}$.

\section{Orthogonality and Jacobi matrices}
\setcounter{equation}{0}

In the sequel $\{ P_n \}_{n \ge 0}$ denotes a sequence of monic
polynomials orthogonal (SMOP)
with respect to a quasi--definite linear functional $u$.

Let $\{ Q_n \}_{n \ge 0}$ be a sequence of monic polynomials with
$\deg Q_n = n$ such that, for $n \ge k+1$,
\begin{equation}\label{Qn}
Q_n(x) = P_n (x)+ a_1 P_{n-1}(x) + \dots  + a_k P_{n-k}(x)
\end{equation}
where the coefficients $\{ a_j \}_{j=1}^k$ are independent of $n$
and $a_k \not= 0$.

Our aim will be to deduce necessary and sufficient conditions in order to
the sequence $\{ Q_n \}_{n \ge 0}$ is orthogonal with respect to a
quasi--definite linear functional $v$ and to give the relation between the
linear functionals $u$ and $v$, via Jacobi matrices.

\begin{proposition} \label{caractincompleta}
Let $\{ P_n \}_{n \ge 0}$ be a sequence of monic orthogonal polynomials with recurrence coefficients
$\{\beta_n \}_{n\ge0}$ and $\{\gamma_n  \}_{n\ge1}$ ($\gamma_n \not= 0$) and let $\{ Q_n \}_{n \ge 0}$ be a
sequence of monic polynomials such that, for $n \ge k+1$,
\begin{equation*}
Q_n (x)= P_n(x)+ a_1 P_{n-1}(x) + \dots  + a_k P_{n-k}(x)
\end{equation*}
where $\{ a_j \}_{j=1}^k$ are constant coefficients and $a_k \not= 0$. Then
$\{ Q_n \}_{n \ge 0}$ is orthogonal with respect to a quasi--definite linear functional if and only
if the following conditions hold

\begin{enumerate}

\item [(i)] For each $j$, $1 \le j \le k$, the polynomials $Q_j$ satisfy a
three term recurrence relation
$xQ_j(x)=Q_{j+1}(x) + \widetilde{\beta_j} Q_j(x) + \widetilde{\gamma_j} Q_{j-1}(x)$, with $\widetilde{\gamma_j} \not= 0$.

\item [(ii)] For $n \ge k+2$
\begin{align*} 
& \gamma_n + a_1 (\beta_{n-1} - \beta_n) = \gamma_{n-k}  \,, \\
& a_{j-1} (\gamma_{n-k} - \gamma_{n-j+1}) = a_j (\beta_{n-j} - \beta_n) \,,
\quad 2\le j \le k
\,.\nonumber
\end{align*}

\item [(iii)] \begin{align*} 
&  \gamma_{k+1} + a_1 (\beta_k - \beta_{k+1}) \not= 0 \\
&a_j \gamma_{k-j+1} + a_{j+1} (\beta_{k-j} -\beta_{k+1}) = a_j^{(k)}[\gamma_{k+1}+ a_1
(\beta_k - \beta_{k+1})] \,, \, 1\le j \le k-1 \,,\nonumber\\
&a_k \gamma_1 = a_k^{(k)}[\gamma_{k+1} + a_1 (\beta_k - \beta_{k+1})] \,,\nonumber
\end{align*}
\end{enumerate}
where $a_j^{(k)}\,, \, j=1, \dots , k \,,$ denotes the coefficient of $P_{k-j}$ in
the Fourier expansion of $Q_k$ in terms of the orthogonal system $\{ P_j 
\}_{j=0}^{k}$.

Moreover, denoting by
$\widetilde{\beta_n}$ and
$\widetilde{\gamma_n}$ the coefficients of the three-term recurrence relation for the 
polynomials $Q_n$ we have for $n \ge k+1$
\begin{equation} \label{widetildes}
 \widetilde{\beta_n} = \beta_n , \quad 
\widetilde{\gamma_n} = \gamma_n + a_1 (\beta_{n-1} - \beta_n)  \,,
\end{equation}
\end{proposition}

\noindent {\bf Proof.}
According to Favard's theorem, the sequence $\{ Q_n \}_{n \ge 0}$ is orthogonal with
respect to a quasi--definite linear functional if and only if, for every $n$,
it satisfies a three--term recurrence relation
$$xQ_n(x)=Q_{n+1}(x) + \widetilde{\beta_n} Q_n (x)+ \widetilde{\gamma_n} Q_{n-1}(x)$$
where $\widetilde{\gamma_n} \not= 0, \, n \ge 1$.
So, condition $(i)$ follows.

Let $n \ge k+2$. From
$\displaystyle x Q_n(x) = x P_n(x) + \sum_{j=1}^k a_j x P_{n-j}(x)$
and a little work involving (\ref{Qn}) and the recurrence relation for the polynomials
$P_n$ it follows that
\begin{align*}
x Q_n (x)&= Q_{n+1}(x) + \beta_n  Q_n (x)+ [\gamma_n + a_1 (\beta_{n-1} - \beta_n)] Q_{n-1}(x)\\
&+ \sum_{j=2}^k \left\{a_j (\beta_{n-j} - \beta_n) - a_{j-1} [\gamma_n - \gamma_{n-j+1} +
a_1(\beta_{n-1} - \beta_n)] \right\} P_{n-j}(x) \\
&- a_k [\gamma_n - \gamma_{n-k} + a_1 (\beta_{n-1} - \beta_n) ] P_{n-(k+1)}(x).
\end{align*}

Then, whenever $n \ge k + 2$, $Q_n$ satisfies a three-term recurrence relation if and only if
the coefficient of $Q_{n-1}$ in the above formula is different from $0$ and the
coefficients of the polynomials $\{ P_j \}_{j= n-(k+1)}^{n-2}$ vanish, i.e., 
\begin{subequations}
\begin{align}
&\gamma_n + a_1 (\beta_{n-1} - \beta_n) \not=0 \label{cns1} \\
&a_{j-1} [\gamma_n - \gamma_{n-j+1} +
a_1(\beta_{n-1} - \beta_n)] = a_j (\beta_{n-j} - \beta_n), \quad j=2,\dots,k  \label{cns2} \\
&\gamma_n + a_1 (\beta_{n-1} - \beta_n) = \gamma_{n-k}  \label{cns3}.
\end{align}
\end{subequations}
Notice that, since $\gamma_n \not=0, n \ge 1,$ (\ref{cns1}) is a consequence of (\ref{cns3}).
Moreover, using (\ref{cns3}), the formula (\ref{cns2}) can be rewritten in the form
\begin{equation*}
a_{j-1} (\gamma_{n-k} - \gamma_{n-j+1}) = a_j (\beta_{n-j} - \beta_n), \quad j=2,\dots,k.
\end{equation*}

So, $(ii)$ holds.

Next, we study the case $n=k+1$. 
Let $Q_k(x) = P_k(x) + \sum_{j=1}^k a_j^{(k)}P_{k-j}(x)$
be the Fourier expansion of $Q_k$ in terms of the orthogonal system $\{ P_n \}_{n \ge
0}$.  Handling in the
same way as above we have

\begin{align*}
&x Q_{k+1}(x) = Q_{k+2}(x) + \beta_{k+1}  Q_{k+1} (x)+ [\gamma_{k+1} + a_1 (\beta_k - \beta_{k+1})] Q_k(x)\\
&+ \sum_{j=1}^{k-1} \left[a_{j+1} (\beta_{k-j} - \beta_{k+1}) - a_j^{(k)}
[\gamma_{k+1} + a_1(\beta_k - \beta_{k+1})] + a_j \gamma_{k-j+1}
\right] P_{k-j}(x),\\
&+[ a_k \gamma_1 - a_k^{(k)}(\gamma_{k+1} + a_1 (\beta_k - \beta_{k+1}))] P_0(x) \,,
\end{align*}
and arguing as in the proof of $(ii)$, $(iii)$ holds.

Finally, (\ref{widetildes}) is an immediate consequence of the precedent results. \ $\Box$

\medskip

\noindent {\bf Remark.} Let us to point out that, because of $(iii)$, the coefficients $\{ a_j^{(k)}
\}_{j=1}^{k}$ are determined by the recurrence parameters
$\{  \beta_n \}_{n\ge0}$ and $\{  \gamma_n  \}_{n\ge1}$ as well as the constants $\{ a_j
\}_{j=1}^{k}$. So, the relation (\ref{Qn}) and the orthogonality of $\{ Q_n \}_{n \ge
k+1}$ fix the polynomial $Q_k$. As a consequence, in the particular case $k=1$, the
sequence $\{ Q_n \}_{n \ge 0}$ is completely determined by (\ref{Qn}) and the
orthogonality property.

\bigskip

Now, we consider two families of monic orthogonal polynomials  $\{
P_n \}_{n\ge 0}$ and $\{ Q_n \}_{n\ge 0}$ with respect to the
quasi-definite linear functionals $u$ and  $v$, respectively,
satisfying the condition (\ref{Qn}). It is well known (see, e.g., \cite{Ma}) that the
relation between the two linear functionals is $u = h_k v$ where $h_k$
is a polynomial of degree $k$.

Writing  ${\bf P} = (P_0, P_1, \dots , P_n, \dots )^T$ and ${\bf Q}
= (Q_0, Q_1, \dots , Q_n, \dots )^T$ for the column vectors associated
with these orthogonal families, and ${\bf J}_P$ and ${\bf J}_Q$ for
the corresponding Jacobi matrices, we get
\begin{equation} \label{matrizJacobi}
x \,\, {\bf P} \, = \, {\bf J}_P \,\, {\bf P} \,, \quad
x \,\, {\bf Q} \, = \, {\bf J}_Q \,\, {\bf Q}   \,.
\end{equation}

If ${\bf M}$ denotes the matrix associated with the change of bases ${\bf Q} = {\bf M} {\bf P}$, then
${\bf M}$ is a lower triangular matrix with  diagonal entries equal to $1$ and zero subdiagonals from the
$(k+1)$--th one.

From (\ref{matrizJacobi}) it follows
${\bf M} \, {\bf J}_P \, {\bf P} \, = x \, {\bf M} \, {\bf P} \, = \, {\bf J}_Q \, {\bf M} \, {\bf P}$
and, therefore,
\begin{equation} \label{MJPJQ}
{\bf M} \, {\bf J}_P \, = \, {\bf J}_Q \, {\bf M} \,.
\end{equation}
From this simple relation it follows straightforward the entries of the matrix ${\bf
J}_Q$.

Moreover, from the equations  (\ref{matrizJacobi}) we get
\begin{equation} \label{submatrizJacobiP}
x ({\bf P})_n \, = \, ({\bf J}_P)_n ({\bf P})_n + P_{n+1} e_{n+1}
\end{equation}
\begin{equation} \label{submatrizJacobiQ}
x ({\bf Q})_n \, = \, ({\bf J}_Q)_n ({\bf Q})_n + Q_{n+1} e_{n+1}
\end{equation}
where $e_{n+1} = (0, \dots , 0, 1)^T \in \RR^{n+1}$. Here, the symbol $({\bf A})_n$
stands for the truncation of any infinite matrix
${\bf A}$ at level $n+1$. Using the relation
(\ref{Qn}), the  representation of the change of bases $ ({\bf Q})_n = ({\bf M})_n \, ({\bf P})_n $
and $(\ref{submatrizJacobiQ})$, we deduce
\begin{equation*}
x ({\bf M})_n ({\bf P})_n
=  ({\bf J}_Q)_n ({\bf M})_n ({\bf P})_n + P_{n+1} e_{n+1} + {\bf L}_n ({\bf P})_n
\end{equation*}
where
\begin{equation*}
{\bf L}_n = \left ( \begin{array}{cccccc}
  0 & \dots  & 0 & 0 & \dots & 0 \\
  . & \dots  & . & . & \dots & . \\
  0 & \dots  & 0 & 0 & \dots & 0  \\
  . & \dots  & . & . & \dots & . \\
  0 & \dots  & 0 & a_k & \dots & a_1 \\
\end{array} \right ) \in \RR^{(n+1,n+1)} \,.
\end{equation*}
Thus,
\begin{equation*}
x ({\bf P})_n = ({\bf M})_n^{-1} \left[ ({\bf J}_Q)_n ({\bf M})_n + {\bf L}_n \right] 
({\bf P})_n + P_{n+1} e_{n+1} \,.
\end{equation*}
Comparing this formula with (\ref{submatrizJacobiP}) we get
\begin{equation*}
({\bf J}_P)_n = ({\bf M})_n^{-1} \left[ ({\bf J}_Q)_n ({\bf M})_n + {\bf L}_n \right]  \,,
\end{equation*}
that is
\begin{equation*}
({\bf J}_Q)_n = ({\bf M})_n \left[ ({\bf J}_P)_n - {\bf L}_n \right] ({\bf M})_n^{-1} \,.
\end{equation*}
This last expression means that $({\bf J}_Q)_n$ is similar to a rank--one perturbation
of the matrix $({\bf J}_P)_n$ and this perturbation is given by the matrix ${\bf L}_n$.
In particular, the zeros of the polynomial $Q_n$ are the zeros of the characteristic 
polynomial of the
matrix $({\bf J}_P)_n - {\bf L}_n$.

\medskip

Next, we are going to describe an explicit algebraic relation
between the Jacobi matrices ${\bf J}_P$ and ${\bf J}_Q$, keeping in
mind basically the relationship between the linear functionals $u$ and
$v$, that is $u=h_k v$.

To do this, we first observe that ${\bf Q} {\bf Q}^T = {\bf M} {\bf P} {\bf P}^T {\bf M}^T \,.$
Writing ${\bf D}_P = \langle u, {\bf P} {\bf P}^T \rangle$ and
${\bf D}_Q = \langle v, {\bf Q} {\bf Q}^T \rangle$ we have
\begin{equation*}
\langle v, h_k {\bf Q} {\bf Q}^T \rangle = \langle h_k v, {\bf Q} {\bf Q}^T \rangle =
\langle u, {\bf Q} {\bf Q}^T \rangle = {\bf M} \langle u, {\bf P} {\bf P}^T \rangle {\bf
M}^T = {\bf M} {\bf D}_P {\bf M}^T.
\end{equation*}
Since $\langle v, h_k {\bf Q} {\bf Q}^T \rangle = \langle v, h_k ({\bf J}_Q) {\bf Q} {\bf
Q}^T \rangle = h_k ({\bf J}_Q) {\bf D}_Q,$ then
\begin{equation} \label{hJQ}
h_k ({\bf J}_Q) = {\bf M} {\bf D}_P {\bf M}^T {\bf D}_Q^{-1}\,.
\end{equation}
On the other hand, from (\ref{MJPJQ}) it follows
\begin{equation} \label{hJQbis}
h_k ({\bf J}_Q) = {\bf M} h_k ({\bf J}_P) {\bf M}^{-1}\,.
\end{equation}
From (\ref{hJQ}) and (\ref{hJQbis}) we deduce
\begin{equation} \label{hJP}
h_k ({\bf J}_P) = {\bf D}_P {\bf M}^T {\bf D}_Q^{-1} {\bf M} \,.
\end{equation}
So, we have a simple algorithm to compute the polynomial $h_k$.
\begin{enumerate}
\item [(1)] From the data ${\bf M}$  and ${\bf J}_P$, we have (\ref{MJPJQ}) and we can 
deduce ${\bf J}_Q \,.$
\item [(2)] From ${\bf J}_P$ and ${\bf J}_Q$ we deduce ${\bf D}_P$ and ${\bf D}_Q ,$ 
respectively.
\item [(3)] Using (\ref{hJP}) and taking into account that $h_k$ is a polynomial of
degree
$k$,
$h_k(x) = c_0 + c_1 x + \dots + c_k x^k$, we get
\begin{equation*}
h_k({\bf J}_P) = c_0 I + c_1 {\bf J}_P + \dots + c_k {\bf J}_P^k
= {\bf D}_P {\bf M}^T {\bf D}_Q^{-1} {\bf M} \,,
\end{equation*}
which is a system of linear equations with $k+1$ unknowns. Notice that the matrices of 
the first
and second terms are $2k+1$ diagonal.
\end{enumerate}
If the  monic polynomials $\{ P_n \}_{n\ge0}$ and $\{ Q_n \}_{n\ge0}$ would be replaced
by the corresponding orthonormal polynomials $\{ \widetilde{P_n} \}_{n\ge0}$ and $\{
\widetilde{Q}_n \}_{n\ge0}$, similar computations would have led to 
\begin{equation*}
h_k({\bf J}_{\widetilde{P}}) =  \widetilde{\bf M}^T \widetilde{\bf M} \,, \qquad
h_k({\bf J}_{\widetilde{Q}}) =  \widetilde{\bf M} \widetilde{\bf M}^T\,,
\end{equation*}
where $\widetilde{\bf M}$ denotes the matrix of the change of bases, that is $\widetilde{\bf Q} =
\widetilde{\bf M} \widetilde{\bf P}$.
This gives us an interesting interpretation of the matrix operation involving the 
linear combination of the orthogonal polynomials
$Q_n(x) = P_n (x)+ a_1 P_{n-1}(x) + \dots  + a_k P_{n-k}(x),\quad n \ge k+1.$

\section{The Case $k=2$}
\setcounter{equation}{0}

Among the classical orthogonal polynomial families, the Chebyshev polynomials are the
unique families such that the sequence of
polynomials $\{Q_n\}_{n\ge0}$ defined by (\ref{Qn})
is orthogonal (see for example \cite {BCG}). But, what happens if the
sequence $\{P_n\}_{n\ge0}$ is not a classical one? 

In this Section, our main goal will
be to describe, for the case $k = 2$, all the families of monic polynomials
$\{P_n\}_{n\ge0}$ orthogonal with respect to a quasi--definite linear functional such that
the new families $\{Q_n\}_{n\ge0}$  are also orthogonal.

\begin{theorem} \label{coeficientesconstantes}
Let $\{ P_n \}_{n\ge0}$ be a SMOP with respect to a quasi--definite linear functional.
Assume that $a_1$ and $a_2$ are real numbers with $a_2\not= 0$ and $ Q_n $ the monic
polynomials defined by
\begin{equation} \label{tresunoconstantes}
Q_n (x)= P_n (x)+ a_1 P_{n-1}(x) + a_2 P_{n-2}(x) \,, \quad n \ge 3 \,.
\end{equation}

Then the orthogonality of the sequence $\{ Q_n \}_{n\ge0}$ depends
on the choice of $a_1$ and $a_2$. More precisely, $\{ Q_n \}_{n\ge0}$ is a SMOP if and only if
$\gamma_3+a_1(\beta_2-\beta_3)\not=0$, and
\item [(i)] if $a_1 = 0$, for $n\ge4$, 
$\beta_n = \beta_{n-2}\quad and \quad \gamma_n = \gamma_{n-2}\, .$

\item [(ii)] if $a_1 \not= 0$ and $a_1^2 = 4a_2$, then for $n\ge2$ \,,
\begin{equation} \label{betanrecurrente}
\beta_n = A + B n + C n^2 \,, \qquad
\gamma_n = D + E n + F n^2 ,
\end{equation}
with $ a_1 C = 2 F\,, \quad a_1 B =2 E - 2 F\,, \quad (A,B,C,D,E,F \in \RR).$

\item [(iii)] if $a_1\not= 0$ and $a_1^2 > 4a_2$, then for $n\ge2$ \,,
\begin{equation*}
\beta_n = A + B {\lambda}^n + C {\lambda}^{-n} \,, \qquad
\gamma_n = D + E {\lambda}^n + F {\lambda}^{-n} \,,
\end{equation*}
with $ a_1 C = (1 + \lambda) F\,, \quad   a_1 \lambda B = (1 + \lambda) E\,, \quad (A,B,C,D,E,F \in \RR) \,,$

\noindent where $\lambda$ is the unique solution in $(-1,1)$ of the equation
$a_1^2
\lambda = a_2 (1+\lambda)^2.$

\item [(iv)] if $a_1\not= 0$ and $a_1^2 < 4a_2$, and let $\lambda = e^{i \theta}$ be the
unique solution  of the equation $a_1^2 \lambda = a_2 (1+\lambda)^2$ with
$\theta \in (0,\pi)$, then for $n\ge2$
\begin{equation*}
\beta_n = A + B e^{in\theta} + \overline B e^{-in\theta} \,, \qquad
\gamma_n = D + E e^{in\theta} + \overline E e^{-in\theta} \,,
\end{equation*}
with $ a_1 \lambda \, B = (1 + \lambda) \,E, \quad (A, D \in \RR\,, B, E \in \CC) \,.$
\end{theorem}

\textbf{Proof.} Applying Proposition \ref{caractincompleta} to the particular case
$k=2$,  we have that $\{ Q_n
\}_{n\ge0}$ is a SMOP if and only if
$\gamma_3+a_1(\beta_2-\beta_3)\not=0$ and, for $n\ge4$,
\begin{equation} \label{unoconstantes}
a_1 ( \gamma_{n-2} - \gamma_{n-1} ) = a_2 (\beta_{n-2} - \beta_n),
\end{equation}
\begin{equation} \label{dosconstantes}
\gamma_n - \gamma_{n-2} = a_1 (\beta_n - \beta_{n-1}).
\end{equation}
Observe that $i)$ follows directly.

In the sequel, we will assume $a_1\not= 0$. From (\ref{unoconstantes}) and 
(\ref{dosconstantes}), we deduce that $\beta_n$ and $\gamma_n$ are
solutions of the difference equation with constant coefficients
\begin{equation} \label{diferec}
y_n + \left(1 - \frac{a_1^2}{a_2}\right) y_{n-1} - \left(1 - \frac{a_1^2}{a_2}\right)
y_{n-2} - y_{n-3} = 0 \,, \quad n\ge5 \,.
\end{equation}
According to the solutions of the associated characteristic equation
\begin{equation} \label{ecuacioncaracteristica}
(\lambda - 1)\left[ \lambda^2 +\left(2 - \frac{a_1^2}{a_2}\right)  \lambda + 1\right] = 0 \,,
\end{equation}
we can analyze three cases (see, for instance, \cite{E}).

(ii) If $a_1^2 = 4a_2$, then $\lambda = 1$ is a root with multiplicity 3 and
therefore
\begin{equation*}
\beta_n = A + B n + C n^2\,, \quad \gamma_n = D + E n + F n^2\,, \quad n \ge 5\,.
\end{equation*}
Note that the obtained expressions for $\beta_n$ and $\gamma_n$ hold also for $n \ge 2,$
just applying (\ref{diferec}) for $n$ equal to $7,$ $6,$ and $5$.

Inserting these expressions of $\beta_n$ and $\gamma_n$ in (\ref{unoconstantes}) and (\ref{dosconstantes})
we have
\begin{equation*}
n [ 2a_1F-a_1^2C] = \frac{1}{2}a_1^2B-a_1E+a_1F, \quad n \ge 4\,,
\end{equation*}
\begin{equation*}
n[4F-2a_1C] = a_1B-a_1C-2 E + 4F \,, \quad n \ge 4\,,
\end{equation*}
which is equivalent to
\begin{equation*}
a_1C-2F=0\,, \qquad   a_1B -2 E + 2F = 0 \,.
\end{equation*}
Moreover, since $\beta_n , \gamma_n \in\RR$, $n\ge1$, it is easy
to check that $A,B,C,D,E,F \in \RR.$

Conversely, the values of $\beta_n$ and $\gamma_n$
given by (\ref{betanrecurrente}), and
the above relations lead, trough (\ref{unoconstantes}) and
(\ref{dosconstantes}), to the orthogonality of the sequence $\{ Q_n
\}.$

\medskip
(iii) and (iv) If $a_1^2 \not= 4a_2$, then
\begin{equation*}
\beta_n = A + B {\lambda}^n + C {\lambda}^{-n}\,, \quad \gamma_n = D + E {\lambda}^n + F {\lambda}^{-n}\,, \quad n \ge 5\,,
\end{equation*}
where $\lambda$ is the unique solution of the equation
(\ref{ecuacioncaracteristica}) such that $\lambda \in (-1,1)$ if
$a_1^2 > 4a_2$ and $\lambda = e^{i\theta}$ with $\theta \in
(0, \pi),$ if $a_1^2 < 4a_2.$ 

Upon applying the same reasoning as in the case (ii) we get that the
previous formulas are also valid for $n \ge 2.$

Inserting these values of $\beta_n$ and $\gamma_n$ in ($\ref{unoconstantes}$) and ($\ref{dosconstantes}$) we
have
\begin{equation*}
\lambda^{2n-2} [ a_1E - a_2B(\lambda+1)]= a_1F\lambda - a_2C(\lambda+1) \,,
\quad n \ge 4\,,
\end{equation*}
\begin{equation*}
\lambda^{2n-2}[a_1B\lambda-(\lambda+1)E] = a_1C-(\lambda+1)F \,, \quad n \ge 4\,.
\end{equation*}
Then, since $\lambda$ is a solution of the equation $a_1^2 \lambda = a_2 (1+\lambda)^2$, we have
that the above  both formulas are equivalent to the following system
\begin{equation*}
 a_1 C = (\lambda+1) F \,, \qquad   a_1\lambda B=(\lambda+1)E \,.
\end{equation*}

Again, the conditions $\beta_n , \gamma_n \in\RR$, $n\ge1$, yield
$A,B,C,D,E,F$ are real numbers in the case (iii) and, in the
case (iv), $A, D \in\RR$ and $B, C, E, F$ are complex numbers with $C=\overline B, \quad
F=\overline E$.$\quad
\Box$

\bigskip
\section{Further remarks and comments}
\setcounter{equation}{0}

After the work of Section 3 it is natural to ask us the following
question: it is possible to give explicitly the SMOP
$\{P_n\}_{n\ge0}$, as well as their orthogonality measure, such
that the sequence $\{Q_n\}_{n\ge0}$ defined by
(\ref{tresunoconstantes}) is also a SMOP? This problem might be
quite hard. In this Section we make some remarks
concerning to it and we show some examples.

\medskip
First, we point out a difference between the cases $k=1$ and $k=2$. Let
$Q_n$ be the monic polynomials defined by
\begin{equation*}
Q_n(x)=P_n(x)+a_1P_{n-1}(x), \quad n\ge2 \,,
\end{equation*}
with $a_1 \not=0$. From Proposition \ref{caractincompleta} written
for $k=1$, it follows that $\{Q_n\}_{n\ge0}$ is a SMOP (see \cite{MP} in a more general
setting), if and only if
\begin{align} \label{gammanbetan}
& \gamma_2+a_1(\beta_1-\beta_2)\not=0 \\
& \gamma_n-\gamma_2=a_1(\beta_n-\beta_2)\,, \quad n\ge3 \,.
\nonumber
\end{align}

Thus, in the case $k=1$, for any sequence of $\{\gamma_n\}_{n\ge1}$
with $\gamma_n\not=0$, if we take $\beta_0, \beta_1 \in \RR$, and $\beta_n$
($n\ge2$) satisfying (\ref{gammanbetan}), we obtain all the SMOP
$\{P_n\}_{n\ge0}$ such that $\{Q_n\}_{n\ge0}$ is also a SMOP. However, in the case
$k=2$, Theorem \ref{coeficientesconstantes} implies that the
recurrence coefficients $\gamma_n$ and $\beta_n$ have to be solutions of the
equation (\ref{diferec}). Therefore, although in both cases we get that $\beta_n$ and $\gamma_n$ have a
similar asymptotic behaviour, roughly speaking, for $k = 2$ there are much less families $\{P_n\}_{n\ge0}$.

\bigskip

{\bf Examples.} 
According to Theorem \ref{coeficientesconstantes},
all the SMOP $\{P_n\}_{n\ge0}$ such that the sequence $\{Q_n \}_{n\ge0}$ where $Q_n = P_n +
a_2 P_{n-2}, n\ge3$ with $a_2\not=0$ is
again a SMOP, satisfy for
$n \ge 4$, $ \beta_n = \beta_{n-2}\quad and \quad \gamma_n = \gamma_{n-2},\,
.$

The families of monic orthogonal polynomials which fulfill these conditions were
explicitly given in terms of Chebyshev polynomials in \cite[Example 2, p. 109]{MP2}.
Observe that this situation corresponds to the case $a_1=0$. However, in the case
$a_1\not=0$, the explicit description of all sequences $\{P_n\}_{n\ge0}$ remains still
open. Besides the four Chebyshev families, we have identified some explicit
solutions, for instance, the continuous big q-Hermite polynomials (see
\cite{KS}).

Whenever $k=1$, an interesting case arises when $\beta_n = \beta_0$, for all $n$ and
$\gamma_n =
\gamma_1, n \ge 2$. In particular, it follows that the only symmetric orthogonal polynomials 
$\{ P_n \}$ such
that the sequence
$P_n + a_1 P_{n-1}$ is also an SMOP are the Chebyshev polynomials (up to a linear change
in the variable).

\end{document}